# Celebrating the first century of ICMI (1908-2008) Some aspects of the history of ICMI

*Ferdinando Arzarello*, *Livia Giacardi*, Department of Mathematics, University of Turin, Italy
*Fulvia Furinghetti*, Department of Mathematics, University of Genoa, Italy
*Marta Menghini*, Department of Mathematics, University of Rome La Sapienza, Italy

**ABSTRACT**

In this paper we report on the events in 2008 that commemorated the Centennial of the International Commission on Mathematical Instruction. This celebration offered the occasion to look back at the history of ICMI and outline the evolution of mathematics education until it achieved its present status as an academic discipline. The years after WWII up to the late 1960s were crucial in this evolution for both the settlement of some institutional aspects (mainly concerning the relationship with mathematicians) and the establishment of new trends of the activities. In this paper we outline – on the basis of unpublished documents - the role of two important figures in those years: Heinrich Behnke and Hans Freudenthal. First as secretary and later as president, Behnke faced the difficult task of reshaping the newborn ICMI after WWII and clarifying the relationship with mathematicians. His mission was completed by Freudenthal, who, as president of ICMI, definitively broke with the past and promoted important initiatives that fostered the emergence of mathematics education as an academic field.

## 1. CELEBRATING THE FIRST CENTURY OF ICMI: A SYMPOSIUM, A BOOK AND A WEBSITE

The IV International Congress of Mathematicians, which took place in Rome from 6 to 11 April 1908, was memorable. An exceptional chronicler, the French mathematician Henri Poincaré, wrote:

> "The number of participants was the highest of any of the preceding Congresses, which is doubtless due to the attraction of the Eternal City, but this is not the only reason. […] France was brilliantly represented […] there were also very distinguished representatives of German science […] No country was absent. […] It goes without saying that Italy had the most and most brilliant representatives […] The sessions were held at the Palazzo Corsini, home of the Accademia dei Lincei […] a beautiful palace in Trastevere […]" (Poincaré, 2008, pp. 19-20, our translation).

It was during this congress that an international commission on the teaching of mathematics was founded; its first president was Felix Klein, an eminent mathematician and promoter of an important reform for the teaching of mathematics in Germany. This commission may be considered the first incarnation of the International Commission on Mathematical Instruction (ICMI).[1]

To celebrate the centennial of the foundation of the ICMI, an international symposium entitled "The First Century of the International Commission on Mathematical Instruction. Reflecting and shaping the world of mathematics education" was held in Rome from 5 to 8 March 2008. Once again, as it did a hundred years ago, Palazzo Corsini, home of the Accademia Nazionale dei Lincei, provided the splendid venue for the congress, along with Palazzo Mattei di Paganica, home of the *Enciclopedia Italiana*.[2]

---

[1] In the first decades of its existence the commission was mainly called *Commission Internationale de l'Enseignement Mathématique* (CIEM) or *Internationale Mathematische Unterrichtskommission* (IMUK). In the following we will use the acronym ICMI to refer to all periods.

[2] The International Programme Committee (IPC), composed of 16 members, was coordinated by Ferdinando Arzarello, while Marta Menghini represented the Organising Committee within the ICP. The permanent website http://www.unige.ch/math/EnsMath/Rome2008/ provides full documentation of the Symposium: the program, the papers presented in the Working Groups, and photos.

The Congress was attended by about 180 participants representing 43 countries. The program included ten plenary lectures, eight parallel lectures, five working groups and a panel discussion. An afternoon was reserved for the Italian teachers, with talks by scholars from Italy and abroad. The talks were broadcast via videoconferences to more than fifty schools throughout Italy; the afternoon session reached more than 1000 teachers.

The last day featured an excursion that recalled that of a century ago, and took the participants to visit the Villa d'Este and Hadrian's Villa in Tivoli, both rich in historical evocations.

The Symposium proceedings have been published by the *Enciclopedia Italiana*, as a volume of their book series *Scienze e Filosofia* (Menghini et al., 2008) and the talks of the Italian afternoon appeared in the journal *Progetto Alice* (n. 25, 2008).

Taking as a point of departure the themes connected to ICMI activities over the course of its hundred year history the symposium sought to identify future directions of research and initiatives for improving mathematics culture in the various countries. The conviction that history is a powerful means not only for giving an account of the past but also for building the future, inspired the activities of the Symposium as well as the publication. The papers in the Proceedings touch on a wide variety of themes: the origins of the ICMI; its rebirth at the end of the 1960s and the emergence of a new field of research; the dialectic between rigour and intuition; the relationships between pure and applied mathematics and the emphasis to be given to each; the interactions between research and practice; the comparison between centres and peripheries of the world; the relationships between mathematics and mathematics teaching; the training of teachers; and the relationship of mathematics education to technology, society and other disciplines. It emerges that ICMI has mirrored the development of mathematics education as a field of study and practice, and stimulated new directions of research, opening new horizons.

The comparison of the "historical" and the "didactical" papers points out the evolution of the research in mathematical education from a collateral aspect of mathematics to an autonomous scientific discipline, whose interactions with mathematics are continuously evolving. The richness of contributions, both in plenary lectures and in working groups, show how varied and deep its current landscape is. The historical papers point out the extent to which ICMI activities in past years had prepared for the growth and development of important themes in didactics of mathematics and the role played by some important figures.

In parallel with the evolution of didactics of mathematics, many papers, and particularly the panel "ICMI's challenges and future", underline a social-geographical evolution from the centre to the peripheries of the world: the rooting of mathematical education in the different cultures gives it a further dimension, which is a richness in today's globalised landscape.

On the occasion of the symposium, a permanent website was created to present the history of the ICMI (see Furinghetti & Giacardi, 2008). Its aim is to delineate the most significant events and the key figures through documents, images and interviews, and make available for the scholars all the tools which are necessary to reconstruct the complex network of relationships that featured this century of history.

The site is divided into six sections: *Timeline; Portrait gallery; Documents; Affiliated Study Groups; International Congresses on Mathematical Education; Interviews and Film clips*. The *Timeline* pinpoints the most important moments in the history of the ICMI: each fact is documented with references to the original sources, in particular to its official organ *L'Enseignement Mathématique* (hereafter EM)[3] with links to its website. Many images, photos and quotations by the protagonists have been inserted. The *Portrait gallery* provides biographic cameos of the ICMI officers during the first hundred years of the Commission who have since passed away, stressing their roles within the ICMI, their publications and their contributions on the problems of teaching. As far as possible the contributions to this section were written by scholars from the country of the officer concerned. Among the *Documents* we find the digitalised versions of the publications of the Central Committee, the texts of the questionnaires proposed for the inquiries and the relative reports, the successive Terms of Reference of ICMI, and the list of the documents kept in the ICMI Archives.[4] The section dedicated to *Affiliated Study Groups* (HPM, ICTMA, IOWME, PME, and WFNMC) presents the history of these groups. The section dedicated to the *International Congresses on Mathematical Education* offers general information about each of them, with bibliographical references, and the Resolutions of the Congress. Finally, *Interviews and Film Clips* provide

---

[3] About this journal and its link with ICMI see (Furinghetti, 2003; 209)

[4] The documents referring to ICMI are in the folders 14 A-G of the IMU files stored at the central Archives of the University of Helsinki. In the following we refer to them as *ICMI Archives* (*IA*). Cf. (Giacardi, 2008b)

the testimony of some of the protagonists of the history of the ICMI – Emma Castelnuovo, Trevor Fletcher, Geoffrey Howson, Maurice Glaymann, Jean-Pierre Kahane, Heinz Kunle, André Revuz and Bryan Thwaites.

2. Glimpses of the first century of ICMI

### 2.1. The five periods of the history of ICMI

The celebration of the hundred years of ICMI[5] has provided an occasion of going through its history, also outlined in (Bass & Hodgson, 2004; Howson, 1984; Lehto, 1998). The book that resulted from the Rome Symposium organized to celebrate the ICMI centenary (Menghini et al., 2008) sheds light on some aspects of the life of this Commission. Donaghue (2008) and Schubring (2008a) document the inception of ICMI, the dissolution and the ephemeral rebirth between the two World Wars. Furinghetti et al. (2008) identify the conditions of ICMI's renaissance after WWII. The recent years, which saw the definitive establishing of mathematics education as an academic discipline, are illustrated in chapters by Bass (2008), Bishop (2008), and Kilpatrick (2008). The website built on the occasion of ICMI's centenary provides detailed information on people and events (Furinghetti & Giacardi, 2008).

It is possible to identify five main periods in the history of ICMI that were produced by both the external events that influenced the life of the Commission as well as by the changing centres of interest and activities of the Commission itself: the early years; the years between the two World Wars; the reconstitution after WWII; the renaissance; the recent decades.

During the early period, from its founding up to WWI, rightly called the "Klein Era", an important international network of national subcommittees was established for the preparation of reports on the state of mathematical instruction as well as on thematic issues. The original aim was to make an inquiry and publish a general report on the current trends in the secondary teaching of mathematics in various countries. However, already in the first meeting after its foundation ICMI acknowledged the need to consider all school levels. The work of the national subcommittees was really impressive. In 1920, at the moment of disbanding, in addition to the eleven Publications of the Central

---

5 Website of ICMI. ttp://www.mathunion.org/ICMI/

Committee there were about 300 reports of the national subcommittees of eighteen countries, for a total number of more than 13,500 pages. At the same time, eight inquiries had been launched and international congresses had been organized. The methodological tenets that underpinned Klein's conception of mathematics teaching and to some extent inspired the early Commission's work concern: bridging the gap between secondary and higher education; the early introduction of the concepts of function and transformation; the applications of mathematics across all the natural sciences; the applications of algebra to geometry and vice-versa; giving importance to the *Approximationsmathematik*, that is, "the exact mathematics of approximate relations"; fostering of intuition in teaching; the approach to topics from a historical perspective favouring a genetic teaching method; the role of elementary mathematics as seen from an advanced standpoint in teacher training.

After WWI the new conditions imposed on official international scientific relations forced international commissions or associations created before the war to either dissolve or reorganize and the shocking decision to ban the researchers of the Central Powers from most international activities was made. ICMI dissolved in 1920 and it was reconstituted only during the ICM in Bologna in 1928, when international collaboration among mathematicians was re-established, reintegrating the countries that had been excluded. However, the Commission was not able to produce new ideas and projects, and was limited to carrying out the old agenda, until WWII forced a second arrest of activities.

In 1952 during the First General Assembly of the reconstituted International Mathematical Union (hereafter IMU), which took place in Rome, the Commission was transformed into a permanent sub-commission of IMU. The new president of ICMI for the period 1952-1954 was Albert Châtelet, dean of the Faculty of Sciences in Paris; the secretary was Heinrich Behnke from the University of Münster, who two years later would become president of ICMI. In the following years, ICMI defined – not without difficulty – several basic structural issues (composition, relationship with the IMU, the organisation of regional international groups, etc.), and established collaborations both scientific and organisational with other associations. These led to a greater internationalism and to the organisation of numerous thematic congresses in various parts of the world, as well as to a broadening of the fields of interest of ICMI.

The actual renaissance and projection into the future took place in the late 1960s. It was Freudenthal who, by establishing almost with a coup de main the tradition of the ICMEs and by founding the new journal *Educational Studies in Mathematics*[6], marked a turning point in the history of the ICMI.

In the last decades there has been an evolution of the relationship between ICMI and IMU, which produced the recent changes in the governance of ICMI, giving more power to the ICMI General Assembly. The focus of the present paper is on some moments that have fostered these changes and, in particular to the two actors in these moments who played a fundamental roles.

## 2.2. Problems and frictions at the rebirth of ICMI

If the first president and enthusiastic promoter was Felix Klein, a considerable role in establishing the Commission was also played by David Eugene Smith, a professor at Teachers College of New York, who was deeply interested in education and in the history of mathematics. Thus the Commission was born of the closest collaboration between mathematicians and educators. The construction of the present status of the Commission and the birth of the discipline "Mathematics Education" are linked to the process of clarification of the relationship with the community of mathematicians and to the carrying out of initiatives independent from this community. The two main characters in the process of constructing new trends were Heinrich Behnke and Hans Freudenthal: both were able to grasp the changes in the world and in mathematics happening in their times and act consequently. Broadening our contribution to the celebration of the centennial we like to reflect on their role stressing the importance of their institutional actions.

Behnke (1898-1979) was a well-known mathematician in the field of complex analysis, professor at the University of Münster, and editor of *Mathematische Annalen* from 1938-1972. Deeply concerned with mathematical education, Behnke invested great energy in teacher training: he had founded the journal *Semesterberichte zur Pflege des Zusammenhangs von Universität und Schule aus den mathematischen Seminaren*, which he edited together with Otto Toeplitz, another mathematician strongly committed to education. The journal was aimed at encouraging the connection between school and university.

In 1955 Behnke also proposed the *Idee völlig utopisch* (Behnke, 1959, p. 148) of realizing of an international encyclopaedia of elementary mathematics, and ardently hoped for the collaboration of mathematics teachers of all levels in order

---

[6] About the foundation of this journal see (Furinghetti, 2008; Hanna & Sidoli, 2002)

to maintain contacts between schools and universities. On that occasion he wrote: "The work done should, if possible, be compiled in a book for the congress of Edinburgh in 1958. The Italian encyclopaedia of elementary mathematics can be regarded as a model in certain ways".[7] In fact in 1958 the German subcommittee of the ICMI produced the first of the five-volume *Grundzüge der Mathematik für Lehrer an Gymnasien sowie für Mathematiker in Industrie und Wirtschaft*. In the preface to the first volume the editors Behnke and Kuno Fladt stress that the work was above all aimed at teachers: "They have always been uppermost in our thoughts. The destiny of future generations of mathematicians depends on their mastery and their love of our science" (p. V-VI). The group of collaborators – which numbered more than 100 members – included scholars not only from Germany but also from Yugoslavia, the Netherlands, Austria and Switzerland, and, significantly, comprised university and high school teachers. Each article has two authors, of which one is a university professor, the other either a high school teacher or someone coming from this career.

As Schubring (2008b) puts it, "[h]aving in so many respects become a true successor of Felix Klein, he eventually followed Klein's footsteps in organizational respect as well". He was the first secretary-general of the renewed Commission, then president of the ICMI from 1955 to 1958, vice-president from 1959 to 1962, and member of the Executive Committee from 1963 to 1970. Behnke's political action was particularly incisive in the period of his secretariat and presidency.[8] That he was completely aware of the problems he had to face in revitalising the Commission emerges from his correspondence: the difficulty of finding mathematicians active in research who were interested in teaching; the difficulty of being recognised in the world of mathematics, and thus how important it was that the work of the commission be visible at the international congresses; the difficulty of obtaining funding; and finally the relevance of the collaboration of mathematics teachers at all levels. In a confidential letter to IMU president Marshall Stone, he wrote:

> "It is a very difficult matter to engage mathematicians, well-known for their research work, into problems of instruction. Most of our colleagues refuse to be active for our commission because they regard this kind of work of little value, and they even neglect to forward circulars. […] The work of our commission reveals its purpose and meaning only when we give lectures and exhibitions at the international mathematical congress".[9]

---

[7] *Program of Work of the International Commission for Mathematical Instruction for the period of 1955/58*, in *IA*, 14 A 1955-1957.

[8] See (Lehto, 1998) and the documents from the ICMI Archives we report below

When Behnke became president he too underlined the need to improve the terms of reference for governing the activities of the Commission and the relationships with IMU; the importance of having regional groups so to decrease the euro-centricity in ICMI; and the need of having ICMI Congresses. In his long report of April 1955 on the activities of ICMI to the IMU president Heinz Hopf, he affirmed:

> "The national sub-commissions suffer from being ruled by university professors, for their influence is predominant through the national adhering organizations … although the number of the university professors in their countries (at least in Europe) represents but a very small part of the teachers of mathematics […] As president of the International Commission of Mathematical Instruction it is my duty to see that the members of the Commission are not university professors only […] the presidency of the national adhering organizations does not appreciate questions of mathematical instruction and inconsiderately uses its national power.
>
> […] The presidency of the IMU has […] to look upon the national sub-commissions – as was the case already before 1914 - as sub-commission of the ICMI, and not of the national adhering organizations. Otherwise the work of the ICMI is made impossible.
>
> […] I regard it a special, honorary mission of the ICMI to establish a contact among the teachers of all levels. The teacher have to get interested in the research work, and those active in the field of research have to get interested in the work of the teachers".[10]

---

[9] Behnke to Stone, Oberwolfach, August 11, 1954, in *IA*, 14A, 1952-1954.
[10] *Report of the president of the International Commission of Mathematical Instruction to the president of the International Mathematical Union*, April 20, 1955, in *IA*, 14A, 1955-1957.

In fact the beginnings of the new commission were not easy and relations with the IMU were characterised by constant friction, derived from the lack of precise terms of reference for governing the activities of the Commission.

Let us quote some passages from the correspondence in ICMI Archives. Stone, IMU president, wrote to Châtelet:

> "It is my understanding that the Commission has proposed an arrangement whereby it will seek the adherence of several nations and set up special national committees in the adhering nations to work with the Commission. I believe that activity of this kind is inappropriate for a Commission of the Union and that it would lead to intolerable confusion as to the relations between the Union, the Commission, and the nations adhering to one or the other".[11]

Stone sent a similar message to Enrico Bompiani, secretary of IMU:

> "There seems to be a great deal of confusion in connection with the ICMI. I hope we can get it cleared up. […] The difficulties […] make me particularly aware of the fact that we need to clarify our procedure for appointing the members of Commission […]".[12]

As well, William Hodge, member of the Executive Committee of IMU wrote to Stone:

> "About ICMI, I agree very strongly that something must be done to curb its activities. At a recent meeting of our national committee very grave concern was expressed at the fact that so many of the Commission's activities were carried on behind our backs and that we were being let in for responsibilities we know nothing about. […] I think it will be necessary to lay down very precise terms of reference for the Commission, and to define its powers very rigidly. It will also be necessary to select a president very carefully".[13]

---

[11] Stone to Châtelet, Chicago, November 3, 1952, in *IA*, 14A, 1952-1954.
[12] Stone to Bompiani, Chicago, July 10, 1953, in *IA*, 14A, 1952-1954.
[13] Stone to Châtelet, Chicago, July 29, 1954, in *IA*, 14A, 1952-1954.

And again Stone to Châtelet:

> "In connection with the Constitution of the National Sub-Commissions, I recall our agreement that each such Sub-Commission is to be in the first place a Sub-Committee of the National Committee for Mathematics in the Country which it represents".[14]

In any case Behnke tried to make up for President Châtelet's lack of initiative, and succeeded in organising the intervention of the ICMI at the International Congress of Mathematicians in Amsterdam in 1954 notwithstanding both the difficulties in relationships between ICMI and IMU, and the resistance on the part of the organising committee of the congress. In a confidential letter to IMU president Stone, he wrote:

> "They [mathematicians] all expressed the idea that lectures on mathematical instruction might not be worthy enough for the Congress. Thus I showed them the reports of previous congresses and pointed out that after 1912 in Cambridge (England) Section VII (history and instruction) was as strongly accentuated as Section II (analysis). […] After a report stating this fact I recommended to rebuilt Section VII. Finally I succeeded […] I was given a highly unfavourable time for the report on the work of our commission, which would never happened in the case of my scientific lectures".[15]

The intention of Behnke was to recreate the climate of fervour and international collaboration that existed during Klein's chairmanship, and he meant "to extend the influence of Section VII (Instruction) so that it will equal the importance it had at the Congress in 1914",[16] the first congress organised by the Commission in Paris.

During the General Assembly of the IMU in The Hague (31 August - 1 September 1954) the Terms of Reference were established, the Executive Committee of the ICMI was renewed and Behnke was nominated president.

---

[14] Hodge to Stone, May 31, 1954, *IMU Archives*, quoted in (Lehto, 1998, p. 111).
[15] Behnke to Stone, Oberwolfach, August 11, 1954, in *IA*, 14A, 1952-1954.
[16] Behnke to Bompiani, Münster, July 8, 1954, in *IA*, 14A, 1952-1954.

According to the new Terms of Reference ICMI had a relatively free hand in its internal organisation, but IMU retained control on important points: the President and the ten members-at-large of ICMI would be elected by the General Assembly of IMU on the nomination of the Union's President. Moreover, the national delegates would be named by each National Adhering Organisation of IMU.

### 2.3. New trends and political issues in 1954-1967

The period from 1954 to 1967 is characterised by important developments both of scientific and organisational kind, which will smooth the way for the ICMI renaissance.

a. New themes to investigate emerged thanks to the interaction with other organisations, such as the Commission Internationale pour l'Étude et l'Amélioration de l'Enseignement des Mathématiques (CIEAEM), which focused on new issues such as: the relevance of psychology in mathematics education; the attention to methodology; the key role of concrete materials; the need to take in consideration all the school levels (from primary to university); empirical research; the relation between mental and mathematical structures (Furinghetti et al., 2008). This different point of view emerges, for example, from the report presented in 1958 at ICM XIII in Edinburgh by Hans Freudenthal (1959) on the comparative study of the methods used in the initiation to geometry: he goes on to a rich and in-depth examination of the teaching subjects and the teaching methods, also considering the impact that psychological and pedagogical research may have on geometrical instruction in the initiating phase.

b. There was an effort to renovate the journal *L'Enseignement Mathématique* by starting a second series with renovated objectives that include space for psychology and methodology in mathematics education:

"This second issue will deal with the subject of mathematics giving special care to the modern theories treating the subject in an easy form studying the methodology and organisation of the teaching, studying the psychological formation of mathematical knowledge, publishing reports on the activity and inquiries of the ICMI. Each issue will contain a bibliographical index".[17]

---

[17] *Note concerning the Review "Enseignement Mathématique"*, in *IA*, 14 A, 1955-1957; see also *EM*, 1955, s. 2, 1, pp. 270-271

c. ICMI tried to gradually reduce its Eurocentric nature by attempting to extend the Commission's activities beyond Europe. In 1955 the Indian Ram Behari was nominated a member of the Executive Committee; in 1956 the ICMI was officially represented by its vice-president Stone in the Conference on Mathematical Instruction in South Asia in Bombay; in 1958 Behnke suggested forming regional groups like the European one to foster international collaboration; in 1961 Stone, president of ICMI, contributed to the organization of the First Inter-American Conference on Mathematics Education in Bogotá.

d. The effort to improve the organization of the Commission led the Executive Committee to discuss the Regulations of the ICMI during the meeting in Brussels on 3 July 1957. In 1958 Behnke proposed a new draft of the by-laws, the main points of which are the following:[18] the reduction to a single representative of each of the National Sub-Commissions; "ICMI is also authorized to accept appropriate organizations as National Sub-Commission even from countries which are not members of IMU"; the right of the National Sub-Commissions to co-opt additional members; "Each National Sub-Commission shall elect a Chairman. Generally, the Chairman shall be the representative to ICMI from his Sub-Commission […] but he also is entitled to delegate a substitute who will have full voting power"; and finally the creation of Regional Groups".

This proposal was criticized by Stone, the new president of the ICMI who wrote to Beno Eckmann, Secretary of the IMU:

> "[…] the way is cleared for the elimination of any real influence in ICMI from the side of the mathematicians who are acquainted with the higher levels of their subject and who are interested in research as well as in teaching and preparation for research".[19]

---

[18] *(Draft) International Commission for Mathematical Instruction (ICMI) New Terms of Reference*, 1958, in *IA*, 14 A, 1958-1960; also *EM*, 1958, s. 2, 4, pp. 216-217.
[19] Stone to Eckmann, Chicago, January 5, 1959, in *IA*, 14 A, 1958-1960.

The New Terms of reference for ICMI would be adopted in 1960: Behnke's proposals were scaled back, and a control on the part of the IMU was sought. Behnke was obviously opposed.[20]

e. Collaborations with other Institutions (OEEC, UNESCO) were sought both for economic reasons and to widen the range of ICMI's influence. During his term as ICMI president Stone (1959-1962) promoted important activities and symposia in collaboration with local organisations, and the Organisation for European Economic Co-operation (OEEC). In 1959, he chaired the influential conference of mathematicians and educators at Royaumont devoted to the new thinking in mathematics and in mathematical education. Other important seminars and symposia were held, from which the guidelines on how to introduce "modern mathematics" into secondary schools emerged.[21] The next president of ICMI, André Lichnérowicz (1963-1966) especially promoted collaboration with UNESCO: UNESCO representatives were officially invited to the congresses organised by the ICMI; ICMI members were consulted as experts by UNESCO; and they were often sent on missions in various countries (Lichnérowicz, 1966). Thanks to the collaboration with UNESCO and other institutions, important international colloquia were organised,[22] and contracts stipulated directly between the ICMI and UNESCO led to the publication of the books in the series *New trends in mathematics teaching*.

## 2.4. The renaissance in the late 1960s and the projection into the future: the role of Freudenthal

In 1967 the presidency of the ICMI passed to Hans Freudenthal (1905-1990). He was a charismatic personality whose broad mathematical knowledge was joined by a profound interest in culture, and who possessed a talent for organisation and the independent spirit necessary to mark a turning point in ICMI activities. In 1963, he entered the Executive Committee of ICMI, as a member until 1966, but from 1967 to 1970 as its president (and thereafter, from 1971 to 1978, as ex-officio member). The story of ICMI was deeply influenced by Freudenthal, who laid the foundations for its renaissance. The inspiring principles of his action were put forward in January 1967 at the UNESCO Colloquium in Lausanne on "Coordination of Instruction of Mathematics and Physics", in which he participated together with other important protagonists of mathematics education of those years, such as Anna Zofia Krigowska and Willy Servais. These principles were published as "Propositions on the teaching of Mathematics" in the first issue of the journal *Educational Studies in mathematics*, which Freudenthal founded in 1968. The main points were:

> "Mathematics constitutes a unique and characteristic activity of human mind. All children have a right to be educated through mathematics".

> "[Mathematics education] must provoke and develop in the first place the capacity of intellectual action instead of merely piling up knowledge".

> "Mathematics develops more and more towards a general science of structures. These structures charge it with a remarkable power of application, information and unification. The knowledge and the mastership [sic] of these structures, its utilization in the grasp of reality are the real objectives of mathematics teaching".

> "The reformation of mathematics teaching has to be considered a permanent process. This implies a continuous retraining of the teachers which is based on regular pedagogical research."[23]

---

[20] Stone to Rolf Nevanlinna, April 5, 1960, in I*A*, 14 A, 1958-1960.

[21] We mention the following meetings: Royaumont, France (23 November- 4 December 1959); Aarhus, Denmark (30 May-2 June 1960); Zagabria – Dubrovnik, Yugoslavia (21 August-19 September 1960); Belgrade, Yugoslavia (19-24 September 1960); Lausanne, Switzerland (26-29 June 1961); Bologna, Italy (4-8 October 1961). See (Furinghetti et al., 2008) and (Giacardi, 2008, 1955-1959, 1960-1966).

[22] We mention the following meetings: Frascati, Italy (8-10 October 1964); Utrecht, Holland (19-23 December 1964); Dakar, Senegal (14-22 January 1965); Echternach, Luxemburg, 30 May - 4 June 1965). See (Giacardi, 2008, 1960-1966).

As Howson (2008, p. 15) observes, "When he assumed the ICMI presidency (January 1967) he was faced with two alternatives: he could carry on as usual

---

[23] *Propositions sur l'enseignement mathématique*, in *IA*, 14B 1967-1974. Propositions on the teaching of mathematics. *Educational Studies in Mathematics*, 1, 1968, 244.

ICME 11 Proceedingsor he could try to break what was by then becoming the ICMI mould. He chose the latter". In fact, from the very beginning Freudenthal intended to open all future ICMI activities to discussion. His request for a permanent ICMI secretariat was rejected by the IMU, which saw neither an urgency nor a purpose for such a move; thus Freudenthal had the greater part of the work of secretary done by his office.[24] He sought and obtained funding beyond that provided by the IMU. He continued the collaboration with UNESCO already well established by his predecessors and stipulated a contract directly with UNESCO for the second volume of *New Trends in Mathematics Teaching* (1970).

The turning point took place in 1967 at the meeting of the Executive Committee of the ICMI held in Utrecht (August 26, 1967), following the colloquium "How to teach mathematics so as to be useful". There, the ideas for the ICME congresses and a new journal were launched. Freudenthal expressed his disappointment about the modality of ICMI participation at the quadrennial Congresses of Mathematicians and supported the idea of a congress of ICMI to be held a year before the ICM. He noted that, in general, the national reports were not useful, so he suggested that future congresses not include the topics of programs and scholastic organisation. He also listed the new subjects for discussion: mathematisation; motivation; how to teach mathematics without a schoolteacher; comparative evaluation of the contents of mathematics courses; criteria of success; evaluation of the results of research in mathematics education; and finally research methodology.

The assembly accepted the project of a Congress of ICMI to be held in 1969. The French delegate Maurice Glaymann proposed holding the congress in France. André Revuz asked for a new journal closer to secondary teachers, because *L'Enseignement Mathématique* was at too high a level. Freudenthal also proposed increasing the number of ICMI members-at-large in order to enliven the Commission, and returned to the question of a permanent secretariat.[25]

On December 2, 1967 IMU secretary Otto Frostman wrote to Freudenthal in an attempt to dissuade him from both initiatives:

---

[24] Frostman to the International Commission on Mathematical Instruction, Djursholm, June 29, 1967, in *IA*, 14B 1967-1974.

[25] EM, 1967, s. 2, 13, pp. 243-246. *Compte-rendu de la séance de la CIEM tenue à Utrecht, le 26 août 1967*, in *IA*, 14B 1967-1974.

80

> "I must admit that I am not too happy about the new pedagogical journal. Do you really think that there is a market for two international journals of that kind (I do not)? If you are not satisfied with *L'Enseignement*, ICMI's official journal, perhaps it would be better to try to reform it. And I am afraid too that in a new journal the "modernizers" of the extreme sort would try to be very busy. At least I ask you to be cautious. I can agree with very much of your criticism of the meetings of ICMI at the International Congresses, but I am not sure that ICMI should isolate itself from those who have, primarily, a scientific interest but who have, nevertheless, very often taken part in the discussions of ICMI. And a special ICMI congress in France in 1969 will cost a lot of money".[26]

On December 20 Freudenthal replied:

> "I would like to reassure you about the new pedagogical journal. The provisional list of editors does not include any "radical". In spite of its name, *Enseignement* has never been a pedagogical journal. Its contributions on education were not pedagogical but organisatory [sic] and administrative. I do not believe it is possible to reorganize a journal so fundamentally".[27]

Frostman wrote to Freudenthal once again on January 2, 1968:

> "I am still a bit afraid that the market will be hard for two publications, even if the new journal will mainly stress other points than *L'Enseignement*."[28]

In March 1969 Frostman complained about not having received any report on ICMI activities[29]. ICMI secretary André Delessert answered him, announcing: the title of the new journal, *Educational Studies in Mathematics*, with Freudenthal as editor and saying that two issues had come out, the first in May 1968, and the second in January 1969; the date and place of the first ICME (Lyon, 24-30 August 1969); the preparation of the journal *Zentralblatt für Didaktik der Mathematik*, published as a collaboration between the ICMI and the *Zentrum für Didaktik der Mathematik* of the University of Karlsruhe. He excuses the delay in providing information by saying that the greater part of the work of secretary is performed by Freudenthal's secretary's office.[30] Thus the IMU was faced with decisions already made. In August 1969 the *First International Congress on Mathematical Education* was held in Lyon. The Congress, attended by 655 active participants from 42 countries, was a big success.[31] The main resolutions concerned: the modernisation of the teaching of mathematics, both in content and method; the collaboration between teachers of mathematics and those of other disciplines; international cooperation; the permanent training of the teachers; the place of "the theory of mathematical education" in universities or research institutes.[32] In the course of the ICMI meeting that took place during the first ICME, Freudenthal explained the reasons for the changes made: although the small congresses dedicated to well-defined topics can be useful, "today we need to go beyond the circle of specialists and reach the teachers, thus large congresses are necessary".[33] (our translation) He further underlined the fact that it is necessary to make it so that all the national sub-committees work and collaborate, and for this it is indispensable that people who are genuinely interested in teaching take part. IMU president Henri Cartan underlined that in any case there had to be retained a section of the ICM dedicated to education and that the members of the ICMI sub-committee had to be designated by the IMU sub-committee.[34]

---

[26] Frostman to Freudenthal, December 2, 1967, in *IA*, 14B 1967-1974.
[27] Freudenthal to Frostman, Utrecht, December 20, 1967, in *IA*, 14B 1967-1974.
[28] Frostman to Freudenthal, January 2, 1968, in *IA*, 14B 1967-1974.
[29] Frostman to Delessert, March 16, 1969, in *IA*, 14B 1967-1974.

On August 1970, during the General Assembly of the IMU in Menton, the IMU President Cartan noted the important work accomplished by outgoing ICMI President Freudenthal, and expressed his desire that the measures that he had begun will come to fruition in the future.

However, he did not even acknowledge the first ICME congress held in Lyon the previous year. During the Assembly James Lighthill was elected ICMI President for the coming four-year term. Shortly before the meeting, Cartan had written to Lighthill suggesting that Freudenthal be kept in the Commission as past president of the ICMI and that a new secretary be chosen who would not be reduced, as Delessert had been, to a "mail box"[35].

---

[30] Delessert to Frostman, Riex, March 22, 1969, in *IA*, 14B 1967-1974.

[31] ICMI Bulletin, 5, 1975, 20-24 and http://www.icmihistory.unito.it/icme1.php.

[32] Cf. ICMI Bulletin, 5, 1975, 20-24 and http://www.icmihistory.unito.it/icme1.php.

[33] *Compte-rendu de la séance de la CIEM tenue à Lyon, le 23 août, à 14 heures, à l'occasion du premier Congrès International de l'Enseignement Mathématique*, in *IA*, 14B 1967-1974.

[34] *Ibidem*.

In the ICMI session held in Nice on the occasion of ICM XVI (1-10 September 1970), the outgoing president Freudenthal presented the decisions of the IMU regarding the composition of the ICMI Executive Committee for the period 1971-1974. From the discussion that followed, two objections emerged: first, all the members at-large of the Commission were appointed by people who were not particularly competent in secondary education (Georges Papy, Freudenthal, Behnke, Đuro Kurepa), second, the members appointed did not represent the various trends in the teaching of elementary mathematics (Papy)[36]. Therefore two important recommendations were formulated: that the regulations which establish the ways that ICMI members are designated had to be modified, and that ICMI members had always to be chosen from among those who are effectively involved with mathematics teaching. Later discussion concerned the organisation of ICME 2: the Congress would take place in Exeter (UK) and would be structured differently from the preceding one, the number of plenary lectures on themes of general interest would be limited, and working groups would be constituted for addressing more specialised topics.

Even at the end of his term, Freudenthal made decisions that irritated the IMU. In fact, even while having to step down as president, he tried to insure that the directions he had opened would be followed with the same aims and guidelines. In October 1970 he sent a letter to the ICMI Executive Committee with a proposal for the Program Committee for ICME 2, which did not include Lighthill, the future president of the ICMI[37].

Three days later Cartan wrote back, letting it be known that it would be the new ICMI who would decide about the organisation of ICME 2; he also requested the rectification of the sentence in the minutes of the session on Sept. 5, 1970, concerning the constitution of the new Executive Committee of the ICMI, underlining that the regulations state that the new Executive Committee had to be designated by the entire new ICMI, that is, after every national sub-committee had designated its representative, and that this, above all, concerned the new president Lighthill[38].

---

[35] Cartan to Lighthill, Die, August 20, 1970, in *IA*, 14B 1967-1974.
[36] *EM*, 1970, s. 2, 16, p. 198.
[37] Freudenthal to the Executive Committee of ICMI, October 11, 1970, in *IA*, 14B 1967-1974.

Freudenthal replied that everything was done in agreement with Lighthill. As to the sentence, in the minutes affirming that ICMI members were often elected by persons who were not effectively competent in mathematics teaching, Freudenthal wrote to Frostman:

> "I admit it looks strong. This, however, reflects the actual discussion in which much stronger terms have been used. The disapproval of the way in which the new members at large were appointed was unanimous. As an attendant to this elections I could only say that the procedure was in complete agreement with the formal regulations. I would suggest that this is taken up as a serious problem by the new Executives of IMU and ICMI".[39]

During his term as president Freudenthal was completely independent with regards to financial matters as well. On November 1970 Frostman wrote to Cartan:

> "I have not paid anything to the ICMI secretariat during the last years [...] in fact I don't have exact information about ICMI's affairs".[40]

In a letter to Lighthill, Cartan wrote that the IMU had provided no funding for ICME in Lyon for the simple reason that nothing was ever requested, and that he had not asked for any funding from UNESCO because Freudenthal had gone to UNESCO directly. He also underlined that the decision to hold ICMI congresses independent of the ICMs was made by Freudenthal without him having ever consulted the IMU, and hoped that Lighthill would establish closer and more confidential relations with the IMU[41].

---

[38] Cartan to Freudenthal, Paris, October 15, 1970, in I*A*, 14B 1967-1974. Writing to Frostman, (Paris, October 15, 1970) Cartan states: "Freudenthal once again worries me […] he overdoes it somewhat by putting the new commission in front of decisions already made" [our translation].
[39] Freudenthal to Cartan, October 19, 1970, and Freudenthal to Frostman, Utrecht, October 23, 1970, in I*A*, 14B 1967-1974.
[40] Frostman to Cartan, November 15, 1970, in I*A*, 14B 1967-1974.
[41] Cartan to Lighthill, Paris, November 30, 1970, in I*A*, 14B 1967-1974.

## 3. AN EPILOGUE

The two important events, the inauguration of the tradition of International Congresses on Mathematical Education (ICMEs), and the launch of journals related to research in Mathematics Education, were made possible thanks to the talent for organisation and the independent spirit of Freudenthal. He realised stable landmarks for its successive development, and the story of ICMI was deeply influenced by him. Continuing the work started by his predecessors in the 1950s he allowed mathematics education to be a discipline in its own right, and not just an appendage to the world of mathematics, so that Steiner (1997, p. 28) was able to write about the ICME-7 in Quebec that

> "for the first time didactics of mathematics showed itself in great clarity as a scientific discipline which under increasing theoretical orientation and empirical foundation is dynamically growing within an international frame of complex cultural, political and interdisciplinary interrelations."

The actions of Behnke and Freudenthal make evident the friction between ICMI and IMU, as well as that between educators and mathematicians active in research, but often inattentive to education, contrasts that have been resolved by the most recent decisions: in fact, starting with the election of the 2010-2012 Executive Committee of the Commission, the election of the ICMI EC is to take place during the General Assembly of ICMI.[42]

Retracing the events that fostered the emancipation of ICMI from the community of mathematicians and the two key figures involved was a further way of celebrating the centenary and expressing our gratitude to them. Bernard of Chartres used to say that we are like dwarfs on the shoulders of giants, so that we can see more than they. So too can we, thanks to our great predecessors, who have contributed to the development of Mathematics Education into an autonomous scientific discipline.

*We are very grateful to all those who have helped us with suggestions and advice, in particular to Michèle Artigue and Bernard Hodgson for their continuous help. We are also grateful to the Associazione Subalpina Mathesis of Torino and to the University of Helsinki for funding the journey to Helsinki to explore ICMI Archives.*

---

[42] See http://www.mathunion.org/organization/ec/procedures-for-election/#ICMI

**APPENDICES**

**1. Excerpt from Report of the president [H. Behnke] of the International Commission of Mathematical Instruction to the president of the International Mathematical Union. April 20, 1955, in IA, 14A, 1955-1957.**

[…] **5.** *The work of the ICMI during 1955/58*
The program of work planned for the ICMI cannot be adopted before the session of the newly constituted Executive Committee has been held. The first session of this Executive Committee will probably take place in Geneva this coming July. But, according to a discussion on the work of the ICMI for the next years at the last session of the former Executive Committee in Paris, Oct. 1954, the following program was suggested:
1. the proposition is to be made to the national sub-commissions to work on the subject of "The Scientific Basis of School mathematics" and to compose for their countries or groups of countries a book for the scientific consultation of the teachers. For this book it is of primary importance that teachers of mathematics of all levels cooperate.

    I regard it is a special, honorary mission of the ICMI to establish a contact among the teachers of all levels. The teacher have to get interested in the research work, and those active in the field of research have to get interested in the work of the teachers. I have already succeeded in being assured of the readiness of cooperation for the second volume of the German ICMI report ("Mathematical Instruction for the early Youth in the Federal Republic of Germany") from professors of the academies for education (Pädagogische Akademien) and through them from the teachers of primary level.

    At the interim meeting in 1956 (symposia for the scientific basis of school mathematics) the experiences shall be compared gathered by the different nations in projecting this book.

    In this context I may mention the suggestion of create an international encyclopedia of elementary mathematics. I do not yet see a way to realize this project because the school systems and therefore the material of instruction deviate too considerably from one another in the different countries. Yet this project will be submitted at the session of the Executive Committee in July. This way it may be possible to approach the suggestion made by M. Stone to create an international work of instruction. This plan

might at first sound simply phantastic for everyone who knows the diversity of national conditions of instruction in different countries. It would, indeed involve an entirely new way of working for the ICMI since, for the first time, it would not simply have to coordinate national work, but would have to realize an important international work.

As a matter of course, considerable financial means would be necessary for the realization of such a project, because the different collaborators would have to be in constant communication during the time of accomplishing this work.

2. The inquiry "The part of Mathematics on Contemporary life" has to be examined more fully and with much more gravity than has been done up to new. The investigation of this bulk of questions is closely connected with the technical development of the different countries.

In America, f. i., there exists a supervision of production at the instant of production. This plays a particularly important part in iron industry of small quantities, thus preventing refuse. The establishment of such a supervision and such a controlling office is, to a high degree, dependent on exact mathematical calculations. Our colleague Ulrich Graf, who died last year, was about to introduce the same establishment in Germany. If this is done on a larger scale in the region of the Ruhr, for example, large numbers of mathematicians will be required. […] Similar questions arise for the use of large-size calculating machines in the industrial field. It is thinkable that, in the coming years, the applications of these machines might expand enormously. This involves the new vocation of the industrial mathematician. The firm Siemens-Halske in Munich has now opened a large department for the development of calculating machines and has already called from Münster four of ours of young doctors of mathematics.

Thus questions are raised which have to be discussed on an international basis.

All pains taken by the ICMI can be summed up by this formula:
To contact people of different qualities and abilities, people of different nations and different teaching professions (as long as they are seriously interested in mathematics) in order to make them work together.
There resides the great obligation and chance for the ICMI.

> I personally try to be an example for this possible, rather comprehensive kind of work

1. with my meetings at Whitsuntide aiming at the maintenance of relation between universities and schools (Pfingsttagungen zur Pflege des Zusammenhangs von Universität and Gymnasien); regular attendance of approximately 250 persons;
2. with the international interim meetings of the ICMI which will be introduced (symposia for school mathematics) and the sessions of the Section at international congresses;
3. with our series of books on mathematical instruction in the different countries;
4. with the national encyclopedias of school mathematics;
   4a. possibly with an international encyclopedia of school mathematics.

## 2. [Memorandum von Herrn Behnke über die Bildung von Gruppen]
in IA, 14A, 1958-1960.

*Suggestions on the subject of forming „Regional Groups" within ICMI*

a. Since its foundation in 1908 in Rome, ICMI's aim has been to compare experiences in the teaching of mathematics in all types of educational establishments, and to discuss possible reforms in the teaching of mathematics. This program includes the following points:

   1. Reports on mathematical instruction
   2. Discussions on eliminating obsolete parts of material hitherto used
   3. Suggestions and discussions on introducing new mathematical points of view into curriculum, (for instance to introduce the concept of "structure" already at Secondary School level)
   4. To establish and cultivate contacts between various types of schools – as far as mathematical teaching is concerned – particularly where pupils progress from one school to the other.

b. In dealing with any particular problem included in this general program, we must be quite clear about the particular age group and the particular educational level of those pupils to whom this problem applies. But this is not easy because conditions vary considerably from country to country, as the national school systems are very often based on different principles. Therefore the work of the *National Sub-Commissions* must be the basis of all life in ICMI. It is then one of the main tasks of ICMI to create and cultivate the exchange of ideas and experiences between the sub-commissions of different countries. This exchange is obviously easiest between those countries where the school systems are similar. Therefore, it is rather natural that the sub-commissions of the countries sharing the old European traditions in educational matters (namely France, Germany, Great Britain, Italy, Scandinavia etc.) – which I shall briefly call the WNE-countries (*W*est and *N*orth *E*uropean) – have up to now found closer contact with one another than with national sub-commissions from other parts of the world. Consequently the activity of ICMI during the 50 years of its existence was mainly concerned with these WNE-countries.

If ICMI now makes the attempt to extend its activities to all parts of the world, then it would be appropriate to form "*groups*" of national sub-commissions, so that countries whose educational systems are similar, belong to the same group. This is in accordance with the resolution passed by the General Assembly of the IMU at St. Andrews, August 1958.[43]

c. At the International Congresses of Mathematicians, ICMI plays a relatively small role, since these congresses are dominated by reports and discussions on matters of research. It might, therefore be more appropriate for ICMI to hold smaller symposia in the years between Congresses. This has, in fact, been the case in the past; especially during the periods 1909-1914 and 1953-1958, such symposia have taken place annually. For financial reasons, it is necessary to restrict each of these meetings to some countries not too far from each other. Thus the geographical aspect must also be considered in the formation of the proposed groups. Fortunately, these two aspects, the similarity of educational systems and the geographical one coincide in the most cases.

Led by these considerations, I make the tentative suggestions that to begin with, the following "Regional Groups" of ICMI must be formed:

---

[43] Cf. (Giacardi 2008a), *1955-1959*.

1. The WNE-countries
2. The East European Countries
3. South East Asia
4. Central and South America
5. Australia and New Zealand.

One might imagine that very large countries, like USA and USSR, would have little interest in joining a regional group.

The organization of the symposia mentioned above would be such that each year a certain group arranges a meeting in which the reports and discussions deal in first place with the particular interests of the members of that group. However representatives of other groups should be invited to take part. With regard to financial arrangements, it remains to be seen whether enough money would be available from IMU or whether the group itself would have to find other resources.

d. The climax of the whole work would be a meeting of all national sub-commissions (at present numbering 23). It would be appropriate for this to take place in connection with each International Congress of Mathematicians.

### 3. Proposition sur l'enseignement mathématique in IA, 14B, 1967-1974.[44]

A la suite de l'enquête de DIALECTICA, les professeurs de mathématique* rassemblés au colloque de Lausanne[45] ont constaté qu'un accord presque général est actuellement réalisé sur les points suivants:

1. La mathématique est une activité inaliénable de l'esprit humain. Tout enfant a le droit d'y être formé.
2. Dans un monde changeant, il convient que cette formation éveille et développe plutôt des aptitudes d'action intellectuelle qu'elle ne fixe des connaissances.

---

[44] These principles were published as "Propositions on the teaching of Mathematics" in the first issue of the journal *Educational Studies in mathematics* (1, 1968, p. 244).

[45] It is the UNESCO Colloquium in Lausanne on "Coordination of Instruction of Mathematics and Physics" (16-10 January 1967).

3. La mathématique évolue de plus en plus vers une science générale des structures. Celles-ci lui confèrent un pouvoir considérable d'application, d'information et d'unification. La connaissance et la maîtrise de ces structures, leur mise en ouvre dans la saisie de la réalité sont les vraies buts de l'enseignement mathématique.
4. Certaines de ces structures ont un caractère élémentaire : il y aurait intérêt à chercher à s'en servir dès l'enfance.
5. Un certain nombre de structures plus élaborées devraient être acquises au terme des études secondaires.
6. La réalisation d'un niveau valable exige une formation mathématique et pédagogique appropriée des maîtres.
7. La réforme de l'enseignement mathématique doit être considérée comme un phénomène permanent. Cela implique une formation continue des maîtres appuyée sur une recherche pédagogique suivie.
8. En ce domaine, une collaboration efficace sur le plan mondial devient indispensable. Il est urgent de fonder un organisme international des informations en matière d'enseignement mathématique.

> Lausanne, le 18 janvier 1967
>
> *A savoir : M. W. Servais (Belgique), M. R. Guy (Canada), M. J. Lichtenberg (Danemark), M. C. Pisot (France), Mme [illegible: P. Gadon ?] et M. A. Renyi (Hongrie), M. C. Cattaneo (Italie), M. H. Freudenthal et M. L. N. H. Bunt (Pays-Bas), M. Z. Krygowska et M.S.Straszewicz (Pologne), M. E. Blanc, M. A. Delessert, M. E. Emery, M. F. Gonseth, M. K. Grimm, M. J. de Siebenthal (Suisse) et M. I. Smolec (Yougoslavie).